\documentclass[12pt]{amsart}
\usepackage{amssymb}
\usepackage{amsthm}
\usepackage{euscript}
\newcommand{\cA}{\EuScript{A}}

\newcommand{\cF}{\EuScript{F}}

\def\ol{\overline}

\theoremstyle{plain}
\newtheorem{lem}[subsection]{Lemma}

\newtheorem{thm}[subsection]{Theorem}

\begin{document}
\title[Burnside Type Theorems]{Burnside type theorems in real and quaternion settings}
\author{Bamdad R. Yahaghi}

\address{Department of Mathematics, Faculty of Sciences, Golestan University, Gorgan 19395-5746, Iran}
\email{bamdad5@hotmail.com, bamdad@bamdadyahaghi.com}
%\thanks{}


\bibliographystyle{plain}

\begin{abstract}

In this note, we consider irreducible semigroups of real,  complex, and  quaternionic matrices with real spectra. We prove Burnside type theorems in the settings of reals and quaternions. First, we prove that  an irreducible semigroup of triangularizable matrices in $M_n(\mathbb{R})$ contains a vector space basis for $M_n(\mathbb{R})$. In other words, $M_n(\mathbb{R})$ is the only irreducible subalgebra of itself  that is spanned by an irreducible semigroup of  triangularizable matrices in $M_n(\mathbb{R})$.   Next, we use this result to show that, up to similarity, $M_n(\mathbb{R})$ is the only irreducible $\mathbb{R}$-algebra in $M_n(\mathbb{H})$  that is spanned by an irreducible semigroup of matrices in $M_n(\mathbb{H})$ with real spectra. Some consequences of our mains results are presented.

\end{abstract}

\maketitle
%\vspace{.5cm}

\bigskip

\begin{section}
{\bf Introduction}
\end{section}

\bigskip

In 1905, W. Burnside proved  a  theorem asserting that a group of \(n \times n\) complex matrices is irreducible iff it contains a vector space basis for $M_n(\mathbb{C})$, equivalently, its linear span is $M_n(\mathbb{C})$, see \cite[Theorem on p. 433]{B}. We extend this result of Burnside to the real field as follows: a semigroup of $n \times n$  triangularizable real matrices is irreducible iff it contains a vector space basis for $M_n(\mathbb{R})$, equivalently, its linear span is $M_n(\mathbb{R})$. In other words, any irreducible semigroup of triangularizable matrices in $M_n(\mathbb{R})$ is absolutely irreducible. Another way to put it is as follows: the algebra $M_n(\mathbb{R})$ is the only irreducible subalgebra of itself  that is spanned by an irreducible semigroup of  triangularizable matrices in $M_n(\mathbb{R})$. Next, we prove the counterpart of our extension of Burnside's Theorem over quaternionic matrices. Finally, we present some consequences of our main results.

Let us begin by establishing  some standard notation and definitions. Let $D$ be a division ring, $F$  a subfield of the center of $D$, denoted by $Z(D)$, and $M_n(D)$  the ring of all $n \times n$ matrices over $D$.  Throughout,  for the most part, the division ring  $D$ will  be the real field, the complex field, or the division ring of quaternions, denoted   by $ \mathbb{R}$, $\mathbb C$, and $\mathbb{H}$, respectively. By $\mathbb F$, we mean $ \mathbb{R}$, or $\mathbb C$, or $\mathbb{H}$. By a semigroup $\mathcal{S} \subseteq M_n(D)$, we mean a set of matrices that is closed under matrix multiplication.
An ideal $\mathcal{J}$ of $\mathcal{S}$ is defined to be a subset of  $\mathcal{S}$ with the property that $SJ \in \mathcal{J}$   and $JS \in \mathcal{J}$ for all $ S \in \mathcal{S}$ and $ J \in \mathcal{J}$.
As is usual, we view the members of $M_n(D)$ as linear transformations acting on the left of $D^n$, where $D^n$ is the right vector space of all $n\times 1$ column vectors.  A family $\cF$ in $M_n(D)$ is called irreducible if the orbit of any nonzero $ x \in D^n$ under $ \mathcal{S}$, the (multiplicative) semigroup generated by  $\cF$, spans $D^n$. When $n > 1$, this is equivalent to the members of $\cF$, viewed as linear transformations on $D^n$, having no common invariant subspace other than the trivial subspaces, namely, $\{0\}$ and $D^n$. We recall that if a semigroup $\mathcal{S}$ of matrices is irreducible, then so is every nonzero semigroup ideal $ \mathcal{J}$ of $\mathcal{S}$ (see \cite[Lemma 2.1.10]{RR1}).  A family $\cF$ in $M_n(D)$ is called reducible if it is not irreducible, or equivalently, $\cF = \{0\}$ or it has a nontrivial invariant
subspace. On the opposite of  irreducibility is triangularizability, when the common invariant subspaces of the members of $\mathcal{S}$ include a maximal subspace chain (of length $n$) in $D^n$, i.e., there are subspaces
$$\{0\} = \mathcal{V}_0 \subseteq  \mathcal{V}_1 \subseteq  \cdots \subseteq \mathcal{V}_n = D^n,$$
where $ \mathcal{V}_j$ is a $j$-dimensional subspace invariant under every $S \in \mathcal{S}$. Such a chain is called a triangularizing chain of subspaces for $\mathcal{S}$. Equivalently, there exists a basis for $D^n$ (called a triangularizing basis for $\mathcal{S}$) relative to which the members of $\mathcal{S}$ have upper triangular matrix representations, or in the language of matrices, there exists an invertible matrix $P \in M_n(D)$ such that $ P^{-1} \mathcal{S} P$ consists of upper triangular matrices, in which case $ P^{-1} \mathcal{S} P$ is called a triangularization of  $\mathcal{S}$ via $P$.
For a  matrix $ A \in M_n(D)$, an element $\lambda \in D$ is said to be a right eigenvalue of $ A \in M_n(D)$ if  there exists a nonzero column vector $ x \in D^n$ such that $ Ax = x \lambda$. By the spectrum of a matrix, we mean the set of all its right eigenvalues.  Note that when $D$ is a field the notions of a right eigenvalue and an eigenvalue of a matrix with entries from $D$ coincide.
In this note,  we are only interested in triangularizable matrices whose spectra are in the center of the division ring. It is easily verified that similar matrices have the same spectra.
It is also easily checked that a matrix $ A \in M_n(D)$ is triangularizable and its spectrum is a subset of the center of $D$ iff the matrix $A$ has a minimal polynomial with coefficients from the center  which splits into linear factors over the center.
By an $F$-algebra $\cA$ in  $ M_n(D)$, we mean a subring of
$ M_n(D)$ that is closed under scalar multiplication by the
elements of the subfield $F$. For a semigroup  $\mathcal S$ in $ M_n(D)$,
we use ${\rm Alg}_F(\mathcal S)$ to denote the $F$-algebra generated by
$\mathcal S$, which is in fact the $F$-linear span of $\mathcal S$. A matrix $A \in M_n(D)$ is called $F$-algebraic if it is algebraic over the subfield $F$.  Let $D, E$ be division rings. The division ring $E$ is called an
extension of $D$ if $D \subset E$ and ${\rm Z}(D) \subset  {\rm
Z}(E)$. A family \(\cF\) of matrices in \(M_n(D)\) is called
 absolutely irreducible if it is irreducible over all
extensions of $D$. If $D=F$, then it follows from Burnside's Theorem that the
family \(\cF\) is absolutely irreducible if and only if \({\rm
Alg}(\cF)= M_n(F)\).  For $ \mathcal{A}, \mathcal{B} \subseteq M_n(D)$, we write $ \mathcal{A} \sim \mathcal{B}$ and say  $\mathcal{A}$ is similar to $\mathcal{B}$ if there exists an invertible matrix $ P \in M_n(D)$ such that $ \mathcal{A} = P^{-1} \mathcal{B} P:= \{P^{-1} BP: B \in \mathcal{B} \} $, in which case we say  that $\mathcal{A}$ is similar to $\mathcal{B}$ via the invertible matrix $ P \in M_n(D)$. It is easily checked that similarity is an equivalence relation on the power set of $M_n(D)$. Naturally, for $A, B \in M_n(D)$, we  write $ A \sim B$ and say  $A$ is (individually) similar to $B$ if  $ \{A\} \sim \{B$\}. Likewise, individual similarity is an equivalence relation on $M_n(D)$.

Throughout, we are mostly dealing  with certain semigroups or $\mathbb R$-algebras in the matrix rings $ M_n(\mathbb{R})$,
 $ M_n(\mathbb{C})$, or  $ M_n(\mathbb{H})$. As is common, we consider the complex field $\mathbb{C}$ and the division ring of quaternions $\mathbb{H}$ as subsets of $M_2(\mathbb{R})$, $M_2(\mathbb{C})$ and $M_4(\mathbb{R})$ via the  embeddings $\varepsilon$, $\varepsilon'$ and $\varepsilon''$, which we call the standard representations of $\mathbb{C}$ and $\mathbb{H}$ in  $M_2(\mathbb{R})$, $M_2(\mathbb{C})$ and $M_4(\mathbb{R})$, respectively.
\begin{eqnarray*}
 \varepsilon (a+ bi) &= & \begin{bmatrix}
a & -b \\
b & a
\end{bmatrix}, \
\varepsilon' (z+ wj) =  \begin{bmatrix}
z & w\\
- \ol{w} & \ol{z}
\end{bmatrix}, \\
\varepsilon'' (a+ bi + cj + dk) &= & \begin{bmatrix}
a & -b & c & -d\\
b & a & d & c \\
-c & -d & a & b\\
d & -c & -b & a
\end{bmatrix}.
 \end{eqnarray*}
We leave it as an exercise to the interested reader to show that the  embeddings $\varepsilon$, $\varepsilon'$ and $\varepsilon''$ are irreducible representations of $\mathbb{C}$ and $\mathbb{H}$ in  $M_2(\mathbb{R})$, $M_2(\mathbb{C})$ and $M_4(\mathbb{R})$, respectively. Likewise, in a natural way, the real algebras $M_n(\mathbb{C})$ and $M_n(\mathbb{H})$ may, respectively, be considered as subsets of $M_{2n}(\mathbb{R})$,  and $M_{2n}(\mathbb{C})$ and $M_{4n}(\mathbb{R})$ via $\varepsilon_n$, $\varepsilon'_n$ and $\varepsilon''_n$,  which are defined as follows
\begin{eqnarray*}
 \varepsilon_n([z_{ij}]_{i, j=1}^n) &=& [  \varepsilon(z_{ij}) ]_{i, j=1}^n, \  z_{ij} \in \mathbb{C}, \\
\varepsilon'_n([q_{ij}]_{i, j=1}^n) &= & [  \varepsilon'(q_{ij}) ]_{i, j=1}^n, \ q_{ij} \in \mathbb{H}, \\
 \varepsilon''_n([q_{ij}]_{i, j=1}^n) &=& [  \varepsilon''(q_{ij}) ]_{i, j=1}^n, \ q_{ij} \in \mathbb{H}.
 \end{eqnarray*}
 Customarily, by writing $\mathbb{C} \subseteq M_2(\mathbb{R})$ and $M_n(\mathbb{H}) \subseteq M_{4n}(\mathbb{R}) $, we actually mean $\varepsilon(\mathbb{C}) \subseteq M_2(\mathbb{R})$ and $M_n(\varepsilon''(\mathbb{H})) \subseteq M_{4n}(\mathbb{R}) $ and so on and so forth.  Perhaps, it is worth mentioning that $M_n(\mathbb{C})$ and $M_n(\mathbb{H})$ are irreducible real subalgebras of $M_{2n}(\mathbb{R})$,  and $M_{2n}(\mathbb{C})$ and $M_{4n}(\mathbb{R})$ (via $\varepsilon_n$, $\varepsilon'_n$ and $\varepsilon''_n$) because so are $\mathbb{C}$ and $\mathbb{H}$ in  $M_2(\mathbb{R})$, $M_2(\mathbb{C})$ and $M_4(\mathbb{R})$ (via $\varepsilon$, $\varepsilon'$ and $\varepsilon''$), respectively; for a detailed proof of a more general case, see \cite[p. 49]{Y3}. If $\mathcal S$ is a semigroup in $M_n(\mathbb{F})$, then the real linear span of $\mathcal S$ forms a real subalgebra of $M_n(\mathbb{F})$, which we call the real algebra spanned by the semigroup $\mathcal S$.

\bigskip

In what follows, we need the following  Wedderburn-Artin type theorem which was proved in \cite{Y2} (see \cite[Theorem 2.2]{Y2}).

\bigskip

%{\bf Theorem 1.1.}
\begin{thm} \label{1.1}
Let $n \in \mathbb{N}$, $D$ be a division ring, $F$ a subfield of its center,
and $\cA$ an irreducible $F$-algebra of $F$-algebraic  matrices in
$ M_n(D)$. Let $r \in \mathbb{N}$  be the smallest nonzero rank present
in $\cA$. Then, the integer $r$ divides $n$ and after a similarity
$\cA = M_{n/r}(\Delta)$, where  $\Delta$ is an irreducible division
$F$-algebra of $F$-algebraic matrices in $ M_r(D)$. In
particular, after a similarity, $\cA =M_n(\Delta_1)$, where $\Delta_1$ is
an $F$-algebraic subdivision ring of $D$, iff $r=1$.
\end{thm}

\bigskip

Also, in order to present our main results, we  need a  characterization of irreducible real subalgebras of $M_n(\mathbb{R})$, $M_n(\mathbb{C})$, and $M_n(\mathbb{H})$. Part (i) of the following theorem is standard and has already appeared in \cite[Theorem 6]{LZ} and \cite[p.415]{LRT}. To the best of our knowledge, part (ii) of the theorem, most likely known by experts, has not appeared anywhere. Part (iii) is a special case of \cite[Theorem 10]{LZ}, which is due to the second named author  and first appeared in \cite{LZ}.  Part (iii) can be thought of as a Burnside type theorem for irreducible $\mathbb R$-algebras of matrices in $ M_n(\mathbb H)$.

\bigskip

%{\bf Theorem 1.2.}
\begin{thm} \label{1.2}
{\rm (i)} Let $n \in \mathbb{N}$, $\cA$ be an irreducible $\mathbb R$-algebra of matrices in
$ M_n(\mathbb R)$, and  $r \in \mathbb{N}$  the smallest nonzero rank present
in $\cA$. Then, the integer $r$, which divides $n$, is $1$, $2$, or $4$, and $\cA = M_{n}(\mathbb R)$, or $\cA \sim M_{n/2}(\mathbb C)$, or
$\cA \sim M_{n/4}(\mathbb H)$, depending on whether $ r= 1$, $r=2$, or $r=4$, respectively.  In
particular, $\cA =M_n(\mathbb R)$ iff $r=1$.

{\rm (ii)}   Let $n \in \mathbb{N}$, $\cA$ be an irreducible $\mathbb R$-algebra of matrices in
$ M_n(\mathbb C)$, and $r \in \mathbb{N}$  the smallest nonzero rank present
in $\cA$. Then, the integer $r$, which divides $n$, is $1$ or $2$, and $\cA = M_{n}(\mathbb C)$ or $\cA \sim M_{n}(\mathbb R)$ if $r=1$ or
$\cA \sim M_{n/2}(\mathbb H)$ if $r=2$.

{\rm (iii)}    Let $n \in \mathbb{N}$, $\cA$ be an irreducible $\mathbb R$-algebra of matrices in
$ M_n(\mathbb H)$, and  $r \in \mathbb{N}$ the smallest nonzero rank present
in $\cA$. Then, $ r=1$ and $\cA = M_{n}(\mathbb H)$, $\cA \sim M_{n}(\mathbb C)$, or
$\cA \sim M_{n}(\mathbb R)$.

\end{thm}

\bigskip

\begin{proof} (i) By Theorem \ref{1.1}, $\mathcal A$ is similar to $M_{n/r}(\Delta)$, where $r \in \mathbb{N}$  is the smallest nonzero rank present
in $\cA$ and $\Delta$ is an irreducible division
$\mathbb R$-algebra of matrices in $ M_r(\mathbb R)$. By Frobenius' Theorem, \cite[Theorem 5.13.12]{La}, $\Delta \cong \mathbb H$ or $\Delta \cong \mathbb C$ or $\Delta \cong \mathbb R$. But $ r = \dim_{\mathbb R} \Delta$ because $\Delta$ is an irreducible division
$\mathbb R$-algebra in $ M_r(\mathbb R)$. Thus, $ r = 1$ or $2$ or $4$. If $r =1$, then $\Delta = \mathbb R$, in which case the assertion follows from Theorem \ref{1.1}. In the other two cases, from the Noether-Skolem Theorem, \cite[p. 39]{D}, we see that   $\Delta \sim \mathbb C  \subseteq M_2(\mathbb R)$ if $r=2$,  or $\Delta \sim \mathbb H \subseteq  M_4(\mathbb R)$ if $ r= 4$. Clearly, the assertion now follows from Theorem \ref{1.1}.

(ii)  By Theorem \ref{1.1}, $\mathcal A$ is similar to $M_{n/r}(\Delta)$, where $r \in \mathbb{N}$  is the smallest nonzero rank present
in $\cA$ and $\Delta$ is an irreducible division
$\mathbb R$-algebra of matrices in $ M_r(\mathbb C)$. By Frobenius' Theorem $\Delta \cong \mathbb H$ or $\Delta \cong \mathbb C$ or $\Delta \cong \mathbb R$. If $\Delta \cong \mathbb R$ or  $\Delta \cong \mathbb C$, from the Noether-Skolem Theorem, \cite[p. 39]{D}, we get that   $\Delta \sim \mathbb R1_{M_r(\mathbb C)} \subseteq M_r(\mathbb C)$  if $\Delta \cong \mathbb R$, or $\Delta \sim \mathbb C 1_{M_r(\mathbb C)} \subseteq M_r(\mathbb C)$  if $\Delta \cong \mathbb C$, implying that $ r=1$ in these two cases because $\Delta$ is irreducible in $ M_r(\mathbb C)$.  If $\Delta \cong \mathbb H$, then  $  \dim_{\mathbb R} \Delta = 4$. On the other hand,  $ \Delta \subseteq  M_r(\mathbb{C})$  is  an irreducible  multiplicative semigroup  in $ M_r(\mathbb{C})$. So, by Burnside's Theorem,  $\Delta$ contains $ r^2$ linearly independent elements over $\mathbb C$. It thus follows  that
$\dim_{\mathbb R} \Delta \geq r^2$, which yields $ r \leq 2$. But $ r > 1$ because $\Delta \cong \mathbb H $ is not commutative. Thus, $ r=2$. Once again, by the Noether-Skolem Theorem, we obtain
$\Delta \sim \mathbb H \subseteq  M_2(\mathbb C)$ if $\Delta \cong \mathbb H$. Therefore, the assertion follows from Theorem \ref{1.1}.

(iii) By Theorem \ref{1.1}, $\mathcal A$ is similar to $M_{n/r}(\Delta)$, where $r \in \mathbb{N}$  is the smallest nonzero rank present
in $\cA$ and $\Delta$ is an irreducible division
$\mathbb R$-algebra of matrices in $ M_r(\mathbb H)$. By Frobenius' Theorem $\Delta \cong \mathbb H$ or $\Delta \cong \mathbb C$ or $\Delta \cong \mathbb R$. This, in view of the Noether-Skolem Theorem, \cite[p. 39]{D}, yields   $\Delta \sim \mathbb H 1_{M_r(\mathbb H)}$ or $\Delta \sim \mathbb C 1_{M_r(\mathbb H)}$ or $\Delta \sim \mathbb R 1_{M_r(\mathbb H)}$, implying that $ r= 1$ because $\Delta$ is irreducible in $M_r(\mathbb H)$. Consequently,
$\Delta \sim \mathbb H$ or $\Delta \sim \mathbb C$ or $\Delta \sim \mathbb R$, finishing the proof by Theorem \ref{1.1}.
\end{proof}

\bigskip

\begin{section}
{\bf Main Results}
\end{section}

\bigskip

We start off with a simple lemma.

\bigskip

\begin{lem} \label{2.1}
Let $\mathcal{S}$ be a semigroup of matrices in ${\rm{M}}_{n}(\mathbb{R})$ and  $\mathcal{S}_{1}=\overline{\mathbb{R} \mathcal{S}}$. Then $\mathcal{S}\subseteq{\rm{M}}_{n}(\mathbb{R})$ is absolutely irreducible  iff $\mathcal{S}_{1}\subseteq{\rm{M}}_{n}(\mathbb{R})$ is.
\end{lem}

\bigskip

\begin{proof}
By Burnside's Theorem,  $\mathcal{S}\subseteq{\rm{M}}_{n}(\mathbb{R})$ is absolutely irreducible semigroup iff ${\rm{Alg}}_{\mathbb{R}}(\mathcal{S})={\rm{M}}_{n}(\mathbb{R})$. Thus,  the assertion is proved as soon as we show that ${\rm{Alg}}_{\mathbb{R}}(\mathcal{S})={\rm{Alg}}_{\mathbb{R}}(\mathcal{S}_{1})$. To this end, as ${\rm{Alg}}_{\mathbb{R}}(\mathcal{S})\subseteq{\rm{Alg}}_{\mathbb{R}}(\mathcal{S}_{1})$, it suffices to show the reverse inclusion. To see this, let $A\in{\rm{Alg}}_{\mathbb{R}}(\mathcal{S}_{1})$ be arbitrary. We get that  $A=\sum_{i=1}^{k}\alpha_{i}S_{i}$, where $k \in \mathbb{N}$, $\alpha_{i}\in\mathbb{R}$, and $S_{i}\in\mathcal{S}_{1}$ ($ 1 \leq i \leq k$). And hence  $S_{i}=\lim_{j}r_{ij}S_{ij}$, where $ r_{ij}\in\mathbb{R} $ and $S_{ij}\in\mathcal{S}$, and $ 1 \leq i \leq k$, $ j \in \mathbb N$. So we can write  $A=\lim_{j}\sum_{i=1}^{k}\alpha_{i}r_{ij}S_{ij}$. But $\sum_{i=1}^{k}\alpha_{i}r_{ij}S_{ij} \in {\rm Alg}_{\mathbb{R}}(\mathcal{S})$ for all $ j \in \mathbb{N}$ and ${\rm Alg}_{\mathbb{R}}(\mathcal{S})$   is a subalgebra, and hence a subspace, of ${\rm{M}}_{n}(\mathbb{R})$, and so it is topologically closed. It thus follows that  $A\in{\rm{Alg}}_{\mathbb{R}}(\mathcal{S})$, proving the reverse inclusion. This completes the proof.
\end{proof}

\bigskip

The following result plays a crucial role in proving our first main theorem.

\bigskip

\begin{thm} \label{2.2}
Let $\mathcal{S}$ be an irreducible  semigroup of triangularizable matrices in ${\rm{M}}_{n}(\mathbb{R})$ that is closed both topologically and under scalar multiplications by  reals, equivalently $\mathcal{S}= \overline{\mathbb{R}\mathcal{S}}$. Then, $\mathcal{S}$ contains a nonzero singular matrix, i.e., a noninvertible matrix other than zero.
\end{thm}

\bigskip

\begin{proof}
By way of contradiction, suppose, on the contrary, that every nonzero member of $\mathcal{S}$ is invertible. We obtain a contradiction by showing that $\mathcal{S} \setminus \{0\}$ is a commutative group of diagonalizable matrices. To this end,  let  $S \in \mathcal{S}\setminus\{0\}$ be arbitrary. Let $S_1=\frac{S}{|\lambda|}$, where $\lambda$ is an eigenvalue of $S$ with the largest absolute value. We show that $S_1$ is diagonalizable and has spectra in $\{-1, 1\}$.
First, note that  $S_1$ is triangularizable, $\sigma(S_1) \subseteq [-1,1]$, and that $\sigma(S_1) \cap \{-1, 1\} \not= \emptyset$, where $\sigma$ stands for the spectrum. Thus $S_1$ is similar to $\begin{bmatrix}
B & 0 \\
0 & C
\end{bmatrix}$ with $\sigma(B) \subseteq \{-1,1\}$ and $\sigma(C)\in (-1,1)\setminus\{0\}$. Next, as $B$ is triangularizable, by the Primary Decomposition Theorem, $B=D+N$, where $D$ is diagonalizable with spectra in $\{-1, 1\}$, $N$ is nilpotent, and $DN=ND$. We prove by contradiction that $N=0$. Suppose otherwise, and let $0 < k <n $ be such that $N^{k}\neq 0$, $N^{k+1}=0$. The binomial expansion yields
$$(D+N)^{n}=D^{n}+\binom{n}{1}D^{n-1}N+\cdots+\binom{n}{k}D^{n-k}N^{k},$$
for all $ n \geq k$. Since $D$ is diagonalizable and has spectra in $\{-1, 1\}$, we see that $D^2= I$, and in particular  $ D^{2j} = I$ for all $j \in \mathbb{N}$.  Setting $ n_j = 2j + k$, we get that
$$\lim_{j}\frac{(D+N)^{n_j}}{\binom{n_j}{k}}=\lim_{j}(D^{n_j-k}N^{k})=N^{k} \not= 0.$$
On the other hand,  if $\rho$ denotes the spectral radius, by Gelfand's formula for the spectral radius, we have $\rho(C)= \lim_n ||C^n||^{1/n}<1$, which, in turn, yields $\lim_{n\to\infty}\|C^{n}\|=0$. Thus,
$$\lim_j\frac{S_1^{n_j}}{\binom{n_j}{k}} \thicksim\begin{bmatrix}
 N^{k} & 0\\
 0 & 0
 \end{bmatrix} \not= 0,$$
 which is impossible because we have obtained  a nonzero and  non-invertible element of $\mathcal S$.
Therefore, $N=0$, and hence $B=D$. Consequently,
$$\lim_j S_1^{2j} \thicksim\begin{bmatrix}
 I & 0\\
 0 & 0
 \end{bmatrix} \not= 0,$$
 implying that  $C$ acts on the zero-dimensional space, i.e., there is no direct summand $C$ in
 $ S_1 \thicksim \begin{bmatrix}
B & 0 \\
0 & C
\end{bmatrix}$.
So, we conclude that $S_1 \thicksim B=D$. Thus, $S_1^{2}\thicksim D^{2}=I$, from which, we obtain $S_1^{2}=I$, and hence $S^2=\lambda ^{2}I$, where $\lambda$ is any eigenvalue of $S$. In other words, we have shown that  $\mathcal{S}_{1}=\lbrace \frac{S}{\rho(S)}:~S\in \mathcal{S}\setminus\{0\}\rbrace$ is a group of involutions. This group, and hence  $\mathcal{S}$, is commutative because $S_1^2 =I$ for all $S_1 \in\mathcal{S}_{1}$. This implies that  $\mathcal{S}$ is simultaneously triangularizable, and hence reducible. This  contradiction completes the proof.
\end{proof}

\bigskip

\noindent {\bf Remark.} In fact, one can prove the following more general result.   {\it Let $\mathcal{S}$ be an irreducible  semigroup of triangularizable matrices in ${\rm{M}}_{n}(\mathbb{F})$ with real spectra  such that $\mathcal{S}= \overline{\mathbb{R} \mathcal{S}}$. Then, $\mathcal{S}$ contains a nonzero singular matrix.} The proof, which is  carried out by mimicking and adjusting the proof above, is omitted for the sake of brevity. It is worth mentioning that individual triangularizability is needed only if  $\mathbb{F}= \mathbb{R}$.

\bigskip

Here is our first main theorem, which is a Burnside type theorem in the setting of reals.

\bigskip

\begin{thm} \label{2.3}
Let $\mathcal{S}$ be an irreducible semigroup of triangularizable matrices in ${\rm{M}}_{n}(\mathbb{R})$. Then ${\rm{Alg}}_{\mathbb{R}}(\mathcal{S})={\rm{M}}_{n}(\mathbb{R})$.
\end{thm}

\bigskip

\begin{proof}
In view of Lemma \ref{2.1}, we may assume, with no loss of generality,  that $\mathcal{S}= \overline{\mathbb{R} \mathcal{S}}$. To prove the assertion, we proceed by induction on $n$, the dimension of the underlying space. If $n=1$, we have nothing to prove. Suppose that the assertion holds for  matrices of size less than $n$. We prove the assertion for matrices of size $n$. By Theorem \ref{2.2}, there is a nonzero singular matrix $S$ in $\mathcal{S}$ so that $ \dim S\mathbb{R}^n < n$.  It is plain that  $S\mathcal{S}|_{S\mathbb{R}^n}$ is an irreducible semigroup of triangularizable linear transformations on the real vector space $S\mathbb{R}^n$. Since $\dim(S\mathbb{R}^n) < n$, we see from the inductive hypothesis that ${\rm Alg}_{\mathbb{R}}(S\mathcal{S}|_{S\mathbb{R}^n})=\mathcal{L}(S\mathbb{R}^n)$. Now, pick a one-rank linear transformation $L$ in ${\rm Alg}_{\mathbb{R}}(S\mathcal{S}|_{S\mathbb{R}^n})$. Consequently, there exists an $A_1 \in {\rm Alg}_{\mathbb{R}}(\mathcal{S})$ such that $ SA_1|_{S\mathbb{R}^n} = L$, from which we obtain $  SA_1S= LS$.
Thus,  $SA_1S$ is a one-rank matrix in ${\rm{Alg}}_{\mathbb{R}}(\mathcal{S})$, and hence  ${\rm{Alg}}_{\mathbb{R}}(\mathcal{S})=M_n(\mathbb{R})$ by Theorem \ref{1.2}(i). This completes the proof.
\end{proof}

\bigskip

Part (ii) of the following theorem is \cite[Theorem 2]{ORR}; it is also a special case of \cite[Theorem 1.1]{Be} with $K = \mathbb{C}$ and $F = \mathbb{R}$.

\bigskip

\begin{thm} \label{2.4}

The following statements are equivalent.

${\rm (i)}$ Let $\mathcal{S}$ be an irreducible semigroup of triangularizable matrices in ${\rm{M}}_{n}(\mathbb{R})$. Then ${\rm{Alg}}_{\mathbb{R}}(\mathcal{S})={\rm{M}}_{n}(\mathbb{R})$.

${\rm (ii)}$  Let $\mathcal{S}\subseteq{\rm{M}}_{n}(\mathbb{C})$ be an irreducible semigroup with spectra in $\mathbb{R}$. Then ${\rm{Alg}}_{\mathbb{R}}(\mathcal{S})\sim{\rm{M}}_{n}(\mathbb{R})$.
\end{thm}

\bigskip

\begin{proof}
``(i) $\Longrightarrow$ (ii)"  Let  $\mathcal{S}\subseteq{\rm{M}}_{n}(\mathbb{C})$ be an irreducible semigroup  with spectra in $\mathbb{R}$. By Theorem \ref{1.2}(ii), we have one of the three cases: (a) ${\rm{Alg}}_{\mathbb{R}}(\mathcal{S})={\rm{M}}_{n}(\mathbb{C})$, or (b) ${\rm{Alg}}_{\mathbb{R}}(\mathcal{S})\sim{\rm{M}}_{n}(\mathbb{R})$, or (c) ${\rm{Alg}}_{\mathbb{R}}(\mathcal{S}) \sim {\rm{M}}_{\frac{n}{2}}(\mathbb{H})$ in which case $2\mid n$. We prove the assertion by refuting the two cases  (a) and (c).  To this end, first suppose (a) holds. But ${\rm{M}}_{n}(\mathbb{C})$ is irreducible in ${\rm{M}}_{2n}(\mathbb{R})$ because so is $\mathbb{C}$ in $M_2(\mathbb{R})$. This in view of (a)  implies that so is ${\rm{Alg}}_{\mathbb{R}}(\mathcal{S})$ in ${\rm{M}}_{2n}(\mathbb{R})$, from which,  we obtain ${\rm{Alg}}_{\mathbb{R}}(\mathcal{S})={\rm{M}}_{2n}(\mathbb{R})$, for (i) holds.  It thus follows that ${\rm{M}}_{n}(\mathbb{C})={\rm{M}}_{2n}(\mathbb{R})$. This is clearly impossible because $\dim_{\mathbb{R}}({\rm{M}}_{n}(\mathbb{C}))\neq\dim_{\mathbb{R}}({\rm{M}}_{2n}(\mathbb{R}))$. So the case (a) is refuted.
Next, suppose (c) holds so that there is a $P\in {\rm{GL}}_{n}(\mathbb{C})$ such that ${\rm{Alg}}_{\mathbb{R}}(\mathcal{S}_{1})={\rm{M}}_{\frac{n}{2}}(\mathbb{H})$, where $\mathcal{S}_{1}=P^{-1}\mathcal{S}P \subseteq {\rm{M}}_{\frac{n}{2}}(\mathbb{H})  \subseteq {\rm{M}}_{2n}(\mathbb{R})$ is an irreducible semigroup with spectra in $\mathbb{R}$. Then again, as $\mathbb{H}$ is irreducible in $M_4(\mathbb{R})$, we see that  ${\rm{M}}_{\frac{n}{2}}(\mathbb{H})$ is irreducible in $ {\rm{M}}_{2n}(\mathbb{R})$. Thus so is $\mathcal{S}_{1}$ in $ {\rm{M}}_{2n}(\mathbb{R})$ because of the equality ${\rm{Alg}}_{\mathbb{R}}(\mathcal{S}_{1})={\rm{M}}_{\frac{n}{2}}(\mathbb{H})$. It thus follows from the hypothesis that ${\rm{Alg}}_{\mathbb{R}}(\mathcal{S}_{1}) ={\rm{M}}_{2n}(\mathbb{R})$, which, in turn, yields ${\rm{M}}_{\frac{n}{2}}(\mathbb{H})={\rm{M}}_{2n}(\mathbb{R})$, which is again clearly impossible. Refuting (a) and (c) establishes (b), which is what we want.

``(ii) $\Longrightarrow$ (i)"   Let $\mathcal{S}$ be an irreducible semigroup of triangularizable matrices in ${\rm{M}}_{n}(\mathbb{R})$ and $r$ be the minimal nonzero rank present in  ${\rm{Alg}}_{\mathbb{R}}(\mathcal{S})$.  By Theorem \ref{1.2}(i), we have one of the three cases (a)   ${\rm{Alg}}_{\mathbb{R}}(\mathcal{S})={\rm{M}}_{n}(\mathbb{R})$, or (b) ${\rm{Alg}}_{\mathbb{R}}(\mathcal{S})\sim {\rm{M}}_{\frac{n}{2}}(\mathbb{C})$, in which case $r=2\mid n$, or (iii)  ${\rm{Alg}}_{\mathbb{R}}(\mathcal{S})\sim{\rm{M}}_{\frac{n}{4}}(\mathbb{H})$, in which case $r= 4\mid n$. Again, we prove the assertion by refuting the two cases (b) and (c). To this end, first suppose (b) holds. Thus,  there is a $P\in {\rm{GL}}_{n}(\mathbb{R})$ such that $P^{-1}{\rm{Alg}}_{\mathbb{R}}(\mathcal{S})P={\rm{M}}_{\frac{n}{2}}(\mathbb{C})$. This yields ${\rm{Alg}}_{\mathbb{R}}(\mathcal{S}_{1})={\rm{M}}_{\frac{n}{2}}(\mathbb{C})$, where $\mathcal{S}_{1}=P^{-1}\mathcal{S}P\subseteq {\rm{M}}_{\frac{n}{2}}(\mathbb{C})$ is an irreducible semigroup with spectra in $\mathbb{R}$. From the hypothesis, we see that  ${\rm{Alg}}_{\mathbb{R}}(\mathcal{S}_{1}) \sim {\rm{M}}_{\frac{n}{2}}(\mathbb{R})$, and hence  ${\rm{M}}_{\frac{n}{2}}(\mathbb{R}) \sim {\rm{M}}_{\frac{n}{2}}(\mathbb{C})$, a  contradiction. This refutes (b). Next, suppose (c) holds.  So there is a $P\in {\rm{GL}}_{n}(\mathbb{R})$ such that ${\rm{Alg}}_{\mathbb{R}}(\mathcal{S}_{1})={\rm{M}}_{\frac{n}{4}}(\mathbb{H})$, where $\mathcal{S}_{1}=P^{-1}\mathcal{S}P\subseteq {\rm{M}}_{\frac{n}{4}}(\mathbb{H})$ is an irreducible semigroup in  ${\rm{M}}_{n}(\mathbb{R})$ with spectra in $\mathbb{R}$.
But  ${\rm{M}}_{\frac{n}{4}}(\mathbb{H})$ is irreducible in ${\rm{M}}_{\frac{n}{2}}(\mathbb{C})$, because so is $\mathbb{H}$ in $M_2(\mathbb{C})$, and hence so is ${\rm{Alg}}_{\mathbb{R}}(\mathcal{S}_{1})$ in ${\rm{M}}_{\frac{n}{2}}(\mathbb{C})$. In other words,  $\mathcal{S}_{1}$ is irreducible in ${\rm{M}}_{\frac{n}{2}}(\mathbb{C})$ and has spectra in  $\mathbb{R}$. This along with the hypothesis implies ${\rm{Alg}}_{\mathbb{R}}(\mathcal{S}_{1}) \sim{\rm{M}}_{\frac{n}{2}}(\mathbb{R})$. Consequently,  ${\rm{M}}_{\frac{n}{4}}(\mathbb{H})\sim{\rm{M}}_{\frac{n}{2}}(\mathbb{R})$. This is a contradiction because $\dim_{\mathbb{R}}({\rm{M}}_{\frac{n}{4}}(\mathbb{H}))\neq\dim_{\mathbb{R}}({\rm{M}}_{\frac{n}{2}}(\mathbb{R}))$. This refutes (c) as well, completing the proof.
\end{proof}

\bigskip

Here is the counterpart of Theorem \ref{2.3} over quaternions, a Burnside type theorem in the setting of quaternions.

\bigskip

\begin{thm}  \label{2.5}
 Let $\mathcal{S}\subseteq{\rm{M}}_{n}(\mathbb{H})$ be an irreducible semigroup of  matrices with spectra in $\mathbb{R}$, equivalently matrices whose minimal polynomials spilt into linear factors over $\mathbb{R}$. Then ${\rm{Alg}}_{\mathbb{R}}(\mathcal{S})\sim{\rm{M}}_{n}(\mathbb{R})$.
\end{thm}

\bigskip

\begin{proof}  Let $\mathcal{S}\subseteq{\rm{M}}_{n}(\mathbb{H})$ be an irreducible semigroup with spectra in $\mathbb R$.  By Theorem \ref{1.2}(iii),  one of the following three cases occurs.  (a)  ${\rm{Alg}}_{\mathbb{R}}(\mathcal{S})\sim{\rm{M}}_{n}(\mathbb{R})$; (b) ${\rm{Alg}}_{\mathbb{R}}(\mathcal{S})\sim{\rm{M}}_{n}(\mathbb{C})$; and (c) ${\rm{Alg}}_{\mathbb{R}}(\mathcal{S})={\rm{M}}_{n}(\mathbb{H})$.
We prove the assertion by rejecting the two cases (b) and (c).
To this end, first suppose (b) holds.
Thus, there is  a $P\in {\rm{GL}}_{n}(\mathbb{H})$ such that $P^{-1}{\rm{Alg}}_{\mathbb{R}}(\mathcal{S})P={\rm{M}}_{n}(\mathbb{C})$.
This yields  ${\rm{Alg}}_{\mathbb{R}}(\mathcal{S}_{1})={\rm{M}}_{n}(\mathbb{C})$, where $\mathcal{S}_{1}=P^{-1}\mathcal{S}P\subseteq {\rm{M}}_{n}(\mathbb{C})$ is an irreducible semigroup of matrices in  ${\rm{M}}_{n}(\mathbb{C})$ (because of the equality above) with spectra in $\mathbb{R}$.
Now, by  Theorems \ref{2.3} and  \ref{2.4}, we get that ${\rm{Alg}}_{\mathbb{R}}(\mathcal{S}_{1}) \sim  {\rm{M}}_{n}(\mathbb{R})$. It thus follows that ${\rm{M}}_{n}(\mathbb{R}) \sim {\rm{M}}_{n}(\mathbb{C})$. This is a contradiction because $\dim_{\mathbb{R}}({\rm{M}}_{n}(\mathbb{R}))\neq\dim_{\mathbb{R}}({\rm{M}}_{n}(\mathbb{C}))$. Next, suppose  (c) holds, i.e., ${\rm{Alg}}_{\mathbb{R}}(\mathcal{S})={\rm{M}}_{n}(\mathbb{H})$. But ${\rm{M}}_{n}(\mathbb{H})$ is irreducible in ${\rm{M}}_{4n}(\mathbb{R})$, because so is $\mathbb{H}$ in $M_4(\mathbb{R})$, and hence so is $\mathcal S$. Consequently, $\mathcal S$ is an irreducible semigroup of triangularizable matrices in ${\rm{M}}_{4n}(\mathbb{R})$.
 Now, by Theorem \ref{2.4}, ${\rm{Alg}}_{\mathbb{R}}(\mathcal{S})={\rm{M}}_{4n}(\mathbb{R})$. It follows that ${\rm{M}}_{n}(\mathbb{H})={\rm{M}}_{4n}(\mathbb{R})$. This is a contradiction because $\dim_{\mathbb{R}}({\rm{M}}_{n}(\mathbb{H}))\neq\dim_{\mathbb{R}}({\rm{M}}_{4n}(\mathbb{R}))$.
\end{proof}

\bigskip

Here is a nice consequence of our main results, namely Theorems \ref{2.3}, \ref{2.4}, and \ref{2.5}. It is worth mentioning that in what follows individual triangularizability is needed only if  $\mathbb{F}= \mathbb{R}$.

\bigskip

\begin{thm}  \label{2.6}
 Let $\mathbb F = \mathbb{R}, \mathbb{C}$, or $\mathbb{H}$, $F$ be a subfield of $\mathbb R$, and  $\mathcal{S}$ an irreducible semigroup of  triangularizable matrices in $M_n(\mathbb{F})$ with spectra in $F$. Then ${\rm{Alg}}_{F}(\mathcal{S})\sim{\rm{M}}_{n}(F)$.
\end{thm}

\bigskip

\begin{proof} In view of Theorems \ref{2.3}, \ref{2.4}, and \ref{2.5}, we may assume that $ \mathbb F = \mathbb{R}$ and that $\mathcal{S}$ is absolutely irreducible. The assertion now follows from \cite[Theorem 1.1]{Be} with $K =\mathbb{R}$ and $F=F$.
\end{proof}

\bigskip

With the above result at our disposal, here is the counterpart of \cite[Theorem B on p. 99]{K} for irreducible semigroups of triangularizable matrices in $M_n(\mathbb{F})$, where $\mathbb F = \mathbb{R}, \mathbb{C}$, or $\mathbb{H}$,  with spectra in a subfield $F$ of $ \mathbb{R}$. The proofs of the next three theorems are given in full details for reader's convenience and for the sake of completeness.

\bigskip

\begin{thm}  \label{2.7}
 Let $ n \in \mathbb{N}$ with $n > 1$,  $\mathbb F = \mathbb{R}, \mathbb{C}$, or $\mathbb{H}$, $F$ be a subfield of $\mathbb R$, and  $\mathcal{S}$ an irreducible semigroup of  triangularizable matrices in $M_n(\mathbb{F})$ with spectra in $F$. Suppose that the set consisting of the sums of the (right) eigenvalues of the members of $\mathcal{S}$ (coutnting multiplicities) has $k$ elements, where $ k \in \mathbb{N}$.  Then $\mathcal{S}$ has  at most $k^{n^2}$ elements, and hence $ k > 1$.
\end{thm}

\bigskip

\begin{proof} In view of Theorem \ref{2.6}, we may assume that $ \mathbb F = F$, that there are $n^2$ linearly independent matrices $S_i$ ($1 \leq i \leq n^2$) in $\mathcal{S}$, and that ${\rm tr}(\mathcal{S})$ has $k$ elements. Since $(S_i)_{i=1}^{n^2}$ forms a linear basis for $M_n(F)$, the homogeneous system
$${\rm tr}(S_i X) = 0 , \  1 \leq i \leq n^2$$
 of $n^2$ linear equations in the $n^2$ unknowns, which are the entries of $X \in M_n(F)$, has only the trivial solution, namely zero. Thus the coefficient matrix of the system is invertible, and hence any nonhomogeneous system of linear equations corresponding to it has a unique solution.  Now, since ${\rm tr}(\mathcal{S})$ has $k$ elements and since every $S \in \mathcal{S}$ is clearly a unique solution of the linear system
$${\rm tr}(S_i S) = y_i , \  1 \leq i \leq n^2, \ y_i \in {\rm tr}(\mathcal{S}),$$
 we see that $\mathcal{S}$ has at most $ k^{n^2}$ elements. This clearly implies $ k > 1$ because $\mathcal{S}$ is irreducible.
\end{proof}

\bigskip

\noindent {\bf Remark.} In the  theorem,  if the set consisting of the sums of the (right) eigenvalues of the members of $\mathcal{S}$ is countable,   then so is $\mathcal{S}$. To see this, imitate the final part of the proof of the following theorem.

\bigskip

In fact with Theorem \ref{2.6} at our disposal, one can state and prove the counterparts of most of the results of \cite{RR2} for irreducible semigroups of triangularizable matrices in $M_n(\mathbb{F})$, where $\mathbb F = \mathbb{R}, \mathbb{C}$, or $\mathbb{H}$,  with spectra in a subfield $F$ of $ \mathbb{R}$, or with spectra in $\mathbb{R}$. For the sake of brevity, we state the counterparts of two main results from \cite{RR2}. The following is the counterpart of \cite[Theorem 1]{RR2}.

\bigskip

\begin{thm}  \label{2.8}
 Let $ n \in \mathbb{N}$ with $n > 1$, $\mathbb F = \mathbb{R}, \mathbb{C}$, or $\mathbb{H}$, $F$ be a subfield of $\mathbb R$, and  $\mathcal{S}$ an irreducible semigroup of  triangularizable matrices in $M_n(\mathbb{F})$ with spectra in $F$ and $\phi$ be an $F$-linear functional on $M_n(\mathbb F)$, which is nonzero on $M_n( F)$. If $\phi(\mathcal{S})$ has $k$ elements, where $ k \in \mathbb{N}$, then $\mathcal{S}$ has  at most $k^{n^2}$ elements, and hence $ k > 1$. Therefore,  $\phi(\mathcal{S})$ is finite iff $\mathcal{S}$ is. Moreover,  $\phi(\mathcal{S})$ is countable iff  $\mathcal{S}$ is.
\end{thm}

\bigskip

\begin{proof} It follows from Theorem \ref{2.6} that there exists an invertible matrix $ P \in M_n(\mathbb{F})$ such that $ M_n(F) = \langle P^{-1} \mathcal{S} P \rangle_F$, where $\langle . \rangle_F$ stands for the $F$-linear span. Let $\phi_F = \phi|_{M_n(F)}$ and $ \mathcal{S}_F = P^{-1} \mathcal{S} P$.  Since $\phi_F$ is a nonzero linear functional on $M_n(F)$,  there exists a nonzero $T_F \in M_n(F)$ such that $ \phi_F (X) = {\rm tr}(T_FX)$ for all $ X \in M_n(F)$. Clearly, $\mathcal{S}_F$ is absolutely irreducible in $M_n(F)$, and hence so is $\mathcal{S}_F T_F \mathcal{S}_F$ (note that  $\mathcal{S}_F T_F \mathcal{S}_F$ is not necessarily is a semigroup ideal of $\mathcal{S}_F$). Consequently, there exists a basis of the form $S_{i1}  T_F S_{i2}$, where  $S_{i1}  , S_{i2} \in \mathcal{S}_F $  and $1 \leq i \leq n^2$, for $M_n(F)$.  Since $(S_{i1}  T_F S_{i2})_{i=1}^{n^2}$ forms a linear basis for $M_n(F)$, the homogeneous system
$${\rm tr}(S_{i1}  T_F S_{i2} X) = 0 , \  1 \leq i \leq n^2$$
 of $n^2$ linear equations in the $n^2$ unknowns, which are the entries of $X \in M_n(F)$, has only the trivial solution, namely zero. Thus, the coefficient matrix of the system is invertible,  and hence any nonhomogeneous system of linear equations corresponding to it has a unique solution.  Now, since ${\rm tr}(T_F \mathcal{S}_F)$ has at most $k$ elements and since every $S \in \mathcal{S}_F$ is clearly a unique solution of the linear system
$${\rm tr}(S_{i1}  T_F S_{i2} S) = y_i , \  1 \leq i \leq n^2, \ y_i= {\rm tr}( T_F S_{i2} S S_{i1} ) \in {\rm tr}(T_F\mathcal{S}_F),$$
 we see that $\mathcal{S}_F$, and hence $\mathcal{S}$ has at most $ k^{n^2}$ elements. This clearly implies $ k > 1$ because $\mathcal{S}$ is irreducible. Finally, if  $\phi(\mathcal{S})$ is countable, then so is  ${\rm tr}(T_F\mathcal{S}_F)$. Now countability of ${\rm tr}(T_F\mathcal{S}_F)$ implies that of $\mathcal{S}_F$, and hence that of $\mathcal{S}$,  because the set of finite subsets of a countable set is countable.
\end{proof}

\bigskip

The following is the counterpart of \cite[Theorem 4]{RR2}.

\bigskip

\begin{thm}  \label{2.9}
 Let $n \in \mathbb{N}$, $\mathbb F = \mathbb{R}, \mathbb{C}$, or $\mathbb{H}$, and  $\mathcal{S}$ an irreducible semigroup of  triangularizable matrices in $M_n(\mathbb{F})$ with spectra in  $\mathbb{R}$ and $\phi$ be an $ \mathbb{R}$-linear functional on $M_n(\mathbb F)$, which is nonzero on $M_n(\mathbb{R})$. If $\phi(\mathcal{S})$  is bounded, then so is $\mathcal{S}$. Therefore, $\mathcal{S}$  is bounded iff $\phi(\mathcal{S})$  is.
\end{thm}

\bigskip

\begin{proof} Just as in the proof of the preceding theorem from Theorem \ref{2.6}, we see that there exists an invertible matrix $ P \in M_n(\mathbb{F})$ such that $ M_n(\mathbb{R}) = \langle P^{-1} \mathcal{S} P \rangle_{\mathbb{R}}$, where $\langle . \rangle_{\mathbb{R}}$ stands for the $\mathbb{R}$-linear span. Let $\phi_{\mathbb{R}} = \phi|_{M_n(\mathbb{R})}$, $ \mathcal{S}_{\mathbb{R}} = P^{-1} \mathcal{S} P$, and $ \mathcal{S}_1 = \overline{[0, 1] \mathcal{S}_{\mathbb{R}} } $. Note that  $ \mathcal{S}_1$ is irreducible and consists of triangularizable matrices. Since $\phi_{\mathbb{R}}$ is a nonzero functional on $M_n(\mathbb{R})$,  there exists a nonzero $T_{\mathbb{R}} \in M_n(\mathbb{R})$ such that $ \phi_{\mathbb{R}} (X) = {\rm tr}(T_{\mathbb{R}}X)$ for all $ X \in M_n(\mathbb{R})$. Clearly, $\phi_{\mathbb{R}}$  is a bounded functional on $ \mathcal{S}_1$. We prove the assertion by showing that $ \mathcal{S}_1$, and hence $\mathcal{S}_{\mathbb{R}}$ and $\mathcal{S}$, are bounded. Suppose by way of contradiction that $ \mathcal{S}_1$ is unbounded so that there exists a sequence $(S_n)_{n=1}^\infty$ in $ \mathcal{S}_1$ with $ \lim_n ||S_n|| = +\infty$. If necessary, by passing to a subsequence, we may assume that $\lim_n \frac{S_n}{||S_n||} = S_0 \in \mathcal{S}_1 \setminus \{0\}$. Then again
$$S S_0 S' = \lim_n \frac{SS_nS'}{||S_n||}, $$
for all $ S, S' \in  \mathcal{S}_1$. This yields
$$\phi_{\mathbb{R}} (SS_0S') = \lim_n   \frac{\phi_{\mathbb{R}}( SS_nS')}{||S_n||}= 0, $$
for all $ S, S' \in  \mathcal{S}_1$. That is the nonzero linear functional $\phi_{\mathbb{R}}$ is zero on the semigroup ideal $\mathcal{S}_1 S_0 \mathcal{S}_1$, which consists of triangularizable matrices. It thus follows from Theorem \ref{2.3} that $ \langle \mathcal{S}_1 S_0 \mathcal{S}_1\rangle_{\mathbb{R}}  = M_n(\mathbb{R}) $. Thus, we get that $\phi_{\mathbb{R}}(\mathcal{S}_1 S_0 \mathcal{S}_1) = 0$, which in turn implies $\phi_{\mathbb{R}}(M_n(\mathbb{R}) )= 0$. This is  a contradiction, completing the proof.
\end{proof}

\bigskip

\noindent {\bf Remark.}  The following example shows that the hypothesis on individual triangularizability of the members of the semigroup having spectra in $\mathbb{R}$ cannot be dropped. Let $n=2$, $\mathbb{F}= \mathbb{R}$ and $\mathcal{S}$ be the multiplicative semigroup generated by $ \mathbb{N} i$, where $i$ is the standard representation of the complex number $i$ as a $2 \times 2$ real matrix. Let $ T = i + {\rm diag}(1, -1)$ and $ \phi: M_n(\mathbb{R}) \longrightarrow  M_n(\mathbb{R})$ be defined by $  \phi(X) = {\rm tr}(TX)$ ($X \in M_n(\mathbb{R})$). It is plain that  $\mathcal{S}$ is an unbounded semigroup in $M_n(\mathbb{R})$ on which the nonzero linear functional $\phi$ is zero, and hence bounded.

\bigskip

%\noindent {\bf Acknowledgement.}


\vspace{2cm}

\begin{thebibliography}{999}


\vspace{1mm}
\bibitem[1]{Be}%{B}
J. Bernik,  The eigenvalue field is a splitting field, {\it Archiv der Mathematik},
Volume 88, Number 6, 2007, 481-490.

\vspace{1mm}
\bibitem[2]{B}
%{\large [B]}
W. Burnside, On the condition of reducibility of any
group of linear substitutions, {\it Proc. London Math. Soc.} 3 (1905),
430-434.



%\vspace{1mm}
%\bibitem[3]{C}
%{\large [DS]}
%P.M. Cohn, {\it Algebra, Volume 3}, 2nd ed., John Wiley \& Sons, Chichester, 1991.

\vspace{1mm}
\bibitem[3]{D}
%{\large [DS]}
P.K. Draxl, {\it Skew Fields}, Cambridge University Press, Cambridge, 1983.



\vspace{1mm}
\bibitem[4]{K}
I. Kaplansky, {\it  Fields and Rings}, 2nd ed., University of Chicago Press, Chicago, 1972.



\vspace{1mm}
\bibitem[5]{La}
%{\large [RR]}
T.Y. Lam, {\it A First Course in Noncommutative Rings}, Springer Verlag, New York, 1991.


\vspace{1mm}
\bibitem[6]{LZ}
T-K. Lee and Y. Zhou, On irreducible and transitive algebras in matrix algebras,
{\it Linear and Multilinear Algebra}, Volume 57, Issue 7, January 2009, pp. 659-672.


%\vspace{1mm}
%\bibitem[8]{L}
%{\large [L]}
%J. Levitzki, {\"U}bber Nilpotente Unterringe,  {\it  Math. Ann.} Vol. {\bf 105}, 1931, 620-627.

\vspace{1mm}
\bibitem[7]{LRT}
V.I. Lomonosov, H. Radjavi, and V.G. Troitsky, Sesquitransitive and localizing operator algebras, {\it Integral Eqns. Operator Theory} 60 (2008), pp. 405–418.

\vspace{1mm}
\bibitem[8]{ORR}
%{\large [ORR]}
M. Omladic, M. Radjabalipour, and H. Radjavi, On semigroups of matrices with traces in a subfield,
{\it Linear Algebra Appl.} {\bf 208/209} (1994), 419-424.


\vspace{1mm}
\bibitem[9]{RR1}
%{\large [RR]}
H. Radjavi and P. Rosenthal, {\it Simultaneous Triangularization}, Springer Verlag, New York, 2000.


\vspace{1mm}
\bibitem[10]{RR2}
%{\large [RR]}
H. Radjavi and P. Rosenthal, Limitations on the size of semigroups of matrices, {\it Semigroup Forum} {\bf 76} (2008), 25-31.

\vspace{1mm}
\bibitem[11]{Rod}
L. Rodman, \emph{ Topics in Quaternion Linear Algebra}, Princeton
University Press, Princeton, New Jersey, 2014.


%\vspace{1mm}
%\bibitem[14]{Ro}
%{\large [RR]}
%L.H. Rowen, {\it Graduate Algebra: Noncommutative View}, American Mathematical Society, Providence, RI, 2008.



\vspace{1mm}
\bibitem[12]{Y1}
%{\large [Y1]}
B.R. Yahaghi, On irreducible semigroups of martices with
traces in a subfield, {\it Linear Algebra Appl.} {\bf 383} (2004),
17-28.


\vspace{1mm}
\bibitem[13]{Y2}
%{\large [Y2]}
B.R. Yahaghi, On $F$-algebras of algebraic matrices over a subfield $F$ of the center of a division ring, {\it Linear Algebra Appl.} {\bf 418} (2006), 599-613.

\vspace{1mm}
\bibitem[14]{Y3}
%{\large [Y3]}
B.R. Yahaghi, {\it Reducibility Results on Operator Semigroups}, Ph.D. Thesis, Dalhousie University, Halifax, Canada, 2002.


\end{thebibliography}
\end{document}